\documentclass[12pt,leqno]{amsart}
\usepackage{amsmath, amssymb}
\usepackage{graphicx,color,hyperref}
\usepackage{verbatim, enumitem}
\newtheorem{theorem}{Theorem}[section]

\newtheorem{question}[theorem]{Question}

\newtheorem{definition}[theorem]{Definition}
\newtheorem{remark}[theorem]{Remark}

\numberwithin{equation}{section}

\def\Aut{\mathop{\rm Aut}\nolimits}

\def\End{\mathop{\rm End}\nolimits}

\def\Ricci{\mathop{\rm Ricci}\nolimits}
\def\exph{\mathop{\rm exph}\nolimits}
\def\rank{\mathop{\rm rank}\nolimits}

\def\tr{\mathop{\rm tr}\nolimits}
\def\PC{\mathop{\rm PC}\nolimits}

\def\bC{{\mathbb C}}
\def\bP{{\mathbb P}}
\def\bR{{\mathbb R}}

\def\cF{{\mathcal F}}
\def\cO{{\mathcal O}}

\begin{document}

\title{Fano manifolds with nef tangent bundles are weakly almost
K\"ahler-Einstein}

\author{Jean-Pierre Demailly}
\date{}


\begin{abstract}
The goal of this short note is to point out that every Fano manifold with a nef tangent bundle possesses an almost K\"ahler-Einstein metric, in a weak sense. The technique relies on a regularization theorem for closed positive $(1,1)$-currents. We also discuss related semistability questions and Chern inequalities.
\end{abstract}

\maketitle

{\hfill {\it dedicated to Professor Ngaiming Mok on the occasion of his sixtieth birthday}}\\

\section{Introduction}
Recall that a holomorphic vector bundle $E$ on a projective manifold $X$ is said to be numerically effective (nef) if the line bundle $\cO_{\bP(E)}(1)$ is nef on $Y=\bP(E)$. Clearly, every homogeneous projective manifold $X$ has its tangent bundle $T_X$ generated by sections, so $T_X$ is nef; morever, by \cite{DPS94}, every compact K\"ahler manifold with $T_X$ nef admits a finite \'etale cover $\tilde X$ that is a locally trivial fibration over its Albanese torus, and the fibers are themselves Fano manifolds $F$ with $T_F$ nef. Hence, the classification problem is essentially reduced to the case when $X$ is a Fano manifold with $T_X$ nef. In this direction, Campana and Peternell \cite{CP91} conjectured in 1991 that Fano manifolds with nef tangent bundles are rational homogeneous manifolds $G/P$, namely quotients of linear algebraic groups by parabolic subgroups. Although this can be checked up to dimension~3 by inspecting the classification of Fano 3-folds by Manin-Iskovskii and Mukai, very little is known in higher dimensions. It would be tempting to use the theory of VMRT's developed by J.M.~Hwang and N.~Mok (\cite{Hwa01}, \cite{Mok08}), since the expected homogeneity property should be reflected in the geometry of rational curves. Even then, the difficulties to be solved remain formidable; see e.g.\ \cite{Mok02} and also \cite{MOSWW} for a recent account of the problem.

On the other hand, every rational homogeneous manifold $X=G/P$ carries a K\"ahler \hbox{metric} that is invariant by a compact real form $G^\bR$ of $G$ (cf.\ \cite{AP86}), and the corresponding Ricci curvature form (i.e.\ the curvature of $-K_{G/P}$) is then a K\"ahler-Einstein metric. A stronger condition than nefness of $T_X$ is the existence of a K\"ahler metric on $X$ whose holomorphic bisectional curvature is nonnegative. N.~Mok \cite{Mok88} characterized those mani\-folds, they are exactly the products of hermitian symmetric spaces of compact type by flat compact complex tori $\bC^q/\Lambda$ and projective spaces $\bP^{n_j}$ with a K\"ahler metric of nonnegative holomorphic bisectional curvature. However, hermitian symmetric spaces of compact type are a smaller class than rational homogeneous manifolds (for instance, they do not include complete flag manifolds), thus one cannot expect these K\"ahler-Einstein metrics to have a nonnegative bisectional curvature in general.

It is nevertheless a natural related question to investigate whether every Fano manifold $X$ with nef tangent bundle actually possesses a K\"ahler-Einstein metric. Our main observation is the following (much) weaker statement, whose proof is based on regularization techniques for closed positive currents \cite{Dem92}, \cite{Dem99}. We denote here by $\PC^{1,1}(X)$ the cone of
positive currents of bidegree $(1,1)$ on $X$ and put $n=\dim_\bC X$.
Though it is a convex set of infinite dimension, the intersection
$\PC^{1,1}(X)\cap c_1(X)$ is compact and metrizable for
the weak topology (see \S$\,$2), and therefore it carries a unique uniform
structure that can be defined by any compatible metric.

\begin{theorem}\label{almost-KE}Let $X$ be a Fano manifold such that the tangent bundle $T_X$ is nef. Then\vskip2pt
\begin{itemize}
\item[\rm (i)] there exists
a family of smoothing operators $(J_\varepsilon)_{\varepsilon\in{}]0,1]}$ that map
every closed positive current $T\in \PC^{1,1}(X)\cap c_1(X)$ to a smooth
closed positive definite $(1,1)$-form
$$
\alpha_\varepsilon=J_\varepsilon(T)\in c_1(X),
$$
in such a way that $J_\varepsilon(T)$ converges weakly to $T$ as $\varepsilon\to 0$ $($uniformly for all~$T\,)$, and
$T\mapsto J_\varepsilon(T)$ is continuous with respect to the weak
topology of currents and the strong topology of $C^\infty$ convergence on
smooth $(1,1)$-forms$\,;$
\item[\rm (ii)] for every $\varepsilon\in{}]0,1]$, there exists
a K\"ahler metric $\omega_\varepsilon$ on $X$ such that
$\Ricci(\omega_\varepsilon)=J_\varepsilon(\omega_\varepsilon)$, in other words
$\omega_\varepsilon $ is ``weakly almost K\"ahler-Einstein''.
\end{itemize}
\end{theorem}

The construction of $J_\varepsilon$ that we have is unfortunately not very
explicit, and we do not even know if $J_\varepsilon$ can be constructed as
a natural linear (say convolution-like) operator. The proof of (ii) is
based on the use of the Schauder fixed point theorem. Let us denote by
$$
0<\rho_{1,\varepsilon}(x)\leq \ldots\leq\rho_{n,\varepsilon}(x)
$$
the eigenvalues of $\Ricci(\omega_\varepsilon)$ with respect to
$\omega_\varepsilon$ at each point $x\in X$, and
$\zeta_{1,\varepsilon},\ldots,\zeta_{n,\varepsilon}\in T_{X,x}$ a corresponding
orthonormal family of eigenvectors with respect to $\omega_\varepsilon$.
We know that
$$
\int_X\sum_{j=1}^n\rho_{j,\varepsilon}\,\omega_\varepsilon^n
=n\int_X\Ricci(\omega_\varepsilon)\wedge\omega_\varepsilon^{n-1}=n\,c_1(X)^n,
\leqno(*)
$$
in particular the left-hand side is bounded. Also, as
$\Ricci(\omega_\varepsilon)-\omega_\varepsilon=J_\varepsilon(\omega_\varepsilon)
-\omega_\varepsilon$ converges weakly to $0$, we know that
$$
\int_X\sum_{j=1}^n(\rho_{j,\varepsilon}-1)\zeta_{1,\varepsilon}^*\wedge
\overline{\zeta_{1,\varepsilon}^*}\wedge u
$$
converges to $0$ for every smooth $(n-1,n-1)$-form $u$ on $X$, hence
the eigenvalues $\rho_{j,\varepsilon}$ ``converge weakly to $1$'' in the sense
that $\sum_{j=1}^n(\rho_{j,\varepsilon}-1)\zeta_{1,\varepsilon}^*\wedge
\overline{\zeta_{1,\varepsilon}^*}$ converges weakly to~$0$ in the space of
smooth $(1,1)$-forms. Therefore,
it does not seem too unrealistic to expect that the family 
$(\omega_\varepsilon)$ is well behaved in the following sense.

\begin{definition}\label{well-behaved}
We say the family $(\omega_\varepsilon)$ of
weak almost K\"ahler-Einstein metrics is {\rm well behaved}
in $L^p$ norm if
$$
\int_X(\rho_{n,\varepsilon}-\rho_{1,\varepsilon})^p\,\omega_\varepsilon^n
\leqno(**)^p
$$
converges to $0$ as $\varepsilon$ tends to $0$.
\end{definition}

\noindent
Standard curvature inequalities then yield the following result.

\begin{theorem}\label{stability}Assume that $X$ is Fano with $T_X$ nef, 
and possesses a family of weak almost K\"ahler-Einstein 
metrics that is well behaved in $L^p$ norm. Then 
\begin{itemize}
\item[\rm (i)] if $p\geq1$, $T_X$ is $c_1(X)$-semistable$\,;$
\smallskip
\item[\rm (ii)] if $p\geq 2$, the
Guggenheimer-Yau-Bogomolov-Miyaoka Chern class inequality
$$
\big[n c_1(X)^2-(2n+2)c_2(X)\big]\cdot c_1(X)^{n-2}\le 0
$$
is satisfied.
\end{itemize}
\end{theorem}

\vskip12pt\noindent
The author is indebted to the referee and to Yury Ustinovskiy for pointing out
some issues in an earlier version of this work. He expresses warm thanks
to Ngaiming Mok, Jun-Muk Hwang, Thomas Peternell, Yum-Tong Siu and
Shing-Tung Yau for inspiring  discussions or viewpoints around these
questions in the last 25 years (or more!).

\section{Proofs}\label{section-main}

\noindent{\sc Existence of regularization operators in}
\ref{almost-KE}$\,$(i). Let us fix
once for all a K\"ahler metric $\beta\in c_1(X)$. It is well known that
$K=\PC^{1,1}(X)\cap c_1(X)$ is compact for the weak topology of currents;
this comes from the fact that currents $T\in K$ have a uniformly bounded mass
$$
\int_X T\wedge\beta^{n-1}=c_1(X)^n.
$$
We also know from \cite{Dem92}, \cite{Dem99} that every such current $T$
admits a family of regularizations $T_\varepsilon\in C^\infty(X)\cap c_1(X)$
converging weakly to $T$, such that $T_\varepsilon\geq -\varepsilon\beta$.
Here, the fact that we can achieve an arbitrary small lower bound uses
in an essential way the assumption that $T_X$ is nef. After replacing
$T_\varepsilon$ by
$T'_\varepsilon=(1-\varepsilon)T_\varepsilon+\varepsilon\beta\geq
\varepsilon^2\beta$, we can assume that $T_\varepsilon$ is a
K\"ahler form. An important issue is that we want to produce a
continuous operator $J_\varepsilon(T)=T_\varepsilon$ defined on
the weakly compact set~$K$. This is clearly the case when the regularization
is produced by convolution, as is done \cite{Dem94}, but in that case one
needs a more
demanding condition on $T_X$, e.g.\ that $T_X$ possesses a hermitian metric with
nonnegative bisectional curvature (i.e.\ that $T_X$ is Griffiths semipositive).
However, the existence of a continuous global operator $J_\varepsilon$ is an
easy consequence of the convexity of~$K$. In fact, the weak topology of $K$
is induced by the $L^2$ topology on the space of potentials
$\varphi$ when writing $T=\beta+dd^c\varphi$ (and taking the quotient by
constants) -- it is also induced by the $L^p$ topology for any $p\in[1,\infty[$,
but $L^2$ has the advantage of being a Hilbert space topology. By compactness,
for every $\delta\in{}]0,1]$ we can then find finitely many currents
$(T_j)_{1\leq j\leq N(\delta)}$ such that the Hilbert balls
$B(T_j,\delta)$ cover $K$. Let $K_\delta$ be the convex hull
of the family $(T_j)_{1\leq j\leq N(\delta)}$ and let
$p_\delta:K\to K_\delta$ be induced by the (nonlinear) Hilbert
projection $L^2/\bR\to K_\delta$.  Since $K_\delta$ is finite
dimensional and is a finite
union of simplices, we can take regularizations $T_{j,\varepsilon}$ of the
vertices $T_j$ and construct a piecewise linear operator
$\tilde J_{\delta,\varepsilon}:K_\delta\to K\cap C^\infty(X)^+$ to
K\"ahler forms, simply by taking linear combinations of the
$T_{j,\varepsilon}$'s. Then
$\smash{\tilde J_{\delta,\varepsilon}}\circ p_\delta(T)$ is 
the continuous regularization operator we need. The uniform weak convergence
to $T$ is guaranteed if we take $J_\varepsilon:=
\smash{\tilde J_{\delta(\varepsilon),\varepsilon}}$ with 
$\varepsilon\ll\delta(\varepsilon)\to 0$, e.g.\ with a step function
$\varepsilon\mapsto \delta(\varepsilon)$ that converges
slowly to $0$ compared to~$\varepsilon$. These operators have the drawback of
being non explicit and a priori non linear.

\begin{question} Is it possible to construct a {\rm linear} regularization
operator $J_\varepsilon$ with the same properties as above, e.g.\ by means
of a convolution process~?
\end{question}

In fact, such a ``linear'' regularization operator $J_\varepsilon$ is constructed
in \cite{Dem94} by putting $T=\beta+dd^c\varphi$ and
$J_\varepsilon(T)=\beta+dd^c\varphi_\varepsilon$ where
$$
\varphi_\varepsilon(z)=\int_{\zeta\in T_{X,z}}
\varphi(\exph_z(\zeta))\,\chi(|\zeta|_h^2/\varepsilon^2)\,dV_h(\zeta)
$$
for a suitable hermitian metrics $h$ on~$T_X$, taking $\exph$ to be the
holomorphic part of the exponential map associated with the Chern
connection~$\nabla_h$. However, the positivity of $J_\varepsilon(T)$ is
in general not preserved, unless one assumes that $(T_X,h)$ is (say)
Griffiths semipositive. As the relation between nefness and Griffiths 
semipositivity is not yet elucidated, one would
perhaps need to extend the above formula to the case
of {\it Finsler metrics}. This is eventually possible by applying some
of the techniques
of~\cite{Dem99} to produce suitable dual Finsler metrics on~$T_X$
(notice that any assumption that $T_X\otimes L$ is ample for some
$\mathbb{Q}$-line bundle $L$ translates into the existence of a strictly
plurisubharmonic Finsler metric on the total space of
$(T_X\otimes L)^*\smallsetminus\{0\}$, so one needs to dualize).\vskip12pt

\noindent{\sc Use of the Schauder fixed point theorem for}
\ref{almost-KE}$\,$(ii). By
Yau's theorem \cite{Yau78}, for every closed $(1,1)$-form
$\rho\in c_1(X)$, there exists a unique K\"ahler metric $\gamma(\rho)\in c_1(X)$
such that $\Ricci(\gamma(\rho))=\rho$. Moreover, by the regularity theory
of nonlinear elliptic operators, the map $\rho\mapsto \gamma(\rho)$ is
continuous in $C^\infty$ topology. We consider the composition
$$
\gamma\circ J_\varepsilon:K\mapsto K\cap C^\infty(X)^+\to
K\cap C^\infty(X)^+\subset K,
\qquad
T\mapsto \rho=J_\varepsilon(T)\mapsto \gamma(\rho).
$$
Since $J_\varepsilon$ is continuous from the weak topology to the strong
$C^\infty$ topology, we infer that $\gamma\circ J_\varepsilon$ is
continuous on $K$ in the weak topology. Now, $K$ is convex and weakly
compact, therefore
$\gamma\circ J_\varepsilon$ must have a fixed point $T=\omega_\varepsilon$
by the theorem of Schauder. By construction $\omega_\varepsilon$ must be
a K\"ahler metric in $c_1(X)$,
since $\gamma\circ J_\varepsilon$ maps $K$ into the space of K\"ahler metrics
contained in $c_1(X)$. This proves Theorem~\ref{almost-KE}.\vskip8pt plus 2pt

\noindent{\sc Semistability of $T_X$} (\ref{stability}$\,$(i)). Let
$\cF\subset \cO(T_X)$ be a coherent subsheaf  such that $\cO(T_X)/\cF$ is
torsion free. We can view $\cF$ as a holomorphic subbundle of $T_X$
outside of an algebraic subset of codimension $2$ in~$X$. The
Chern curvature tensor $\Theta_{T_X,\omega_\varepsilon}= \frac{i}{2\pi}\nabla^2$
satisfies the Hermite-Einstein condition
$$
\Theta_{T_X,\omega_\varepsilon}\wedge \omega_\varepsilon^{n-1}=
\frac{1}{n}\rho_\varepsilon\,\omega_\varepsilon^n,
$$
where $\rho_\varepsilon\in C^\infty(X,\End(T_X))$ is the Ricci operator
(with the same eigenvalues $\rho_{j,\varepsilon}$ as $\Ricci(\omega_\varepsilon)$).
It is well known that the curvature of a subbundle is always bounded above by the restriction of the full curvature tensor, i.e.\
$\Theta_{\cF,\omega_\varepsilon}\leq\Theta_{T_X,\omega_\varepsilon|\cF}$ (say, in the
sense of Griffiths positivity, viewing the curvature tensors as hermitian
forms on $T_X\otimes\cF$). By taking the trace with respect to
$\omega_\varepsilon$, we get
$$
\Theta_{\cF,\omega_\varepsilon}\wedge \omega_\varepsilon^{n-1}\leq
\frac{1}{n}\rho_{\varepsilon|\cF}\,\omega_\varepsilon^n.
$$
By the minimax principle, if $r=\rank\cF$, the eigenvalues of
$\rho_{\varepsilon|\cF}$ are bounded above by
$$
\rho_{n-r+1,\varepsilon}\leq\ldots\leq\rho_{n,\varepsilon},
$$
hence
$$
\int_X c_1(\cF)\wedge \omega_\varepsilon^{n-1} =
\int_X \tr_{\cF}\Theta_{\cF,\omega_\varepsilon}\wedge\omega_\varepsilon^{n-1}
\leq \frac{1}{n}\int_X
\sum_{j=1}^r\rho_{n-r+j,\varepsilon}\,\omega_\varepsilon^n.
$$
The left hand side is unchanged if we replace $\omega_\varepsilon$ by
$\beta\in c_1(X)$, and our assumption $(**)^p$ for $p=1$
(cf.~Definition~\ref{well-behaved}) implies
$$
\int_X c_1(\cF)\wedge \beta^{n-1} \leq\lim_{\varepsilon\to 0}
\frac{1}{n}\int_X
\sum_{j=1}^r\rho_{n-r+j,\varepsilon}\,\omega_\varepsilon^n
=\lim_{\varepsilon\to 0}\frac{r}{n^2}\int_X
\sum_{j=1}^n\rho_{j,\varepsilon}\,\omega_\varepsilon^n=\frac{r}{n}
\int_Xc_1(X)\wedge\beta^{n-1}
$$
by $(*)$. This means that $T_X$ is $c_1(X)$-semistable.

\vskip8pt plus 2pt

\noindent{\sc Chern class inequality} (\ref{stability}$\,$(ii)).
Let us write
$$
\Theta_{T_X,\omega_\varepsilon}
=(\theta_{\alpha\beta})_{1\leq \alpha,\beta\leq n}
=\Big(\sum_{j,k}\theta_{\alpha\beta jk}dz_j\wedge d\overline z_k
\Big)_{1\leq \alpha,\beta\leq n}
$$
as an $n\times n$
matrix of $(1,1)$-forms with respect to an orthonormal frame
of $T_X$ that diagonalizes the Ricci operator $\rho_\varepsilon$
(with respect to the metric $\omega_\varepsilon$). A standard calculation
yields
\begin{eqnarray*}
{}\kern45pt&&{}\kern-90pt 2\Big[n\,c_1(X)_h^2{}-(2n+2)\,c_2(X)\Big]\wedge
\frac{\omega_\varepsilon^{n-2}}{(n-2)!}\\
&&\kern-26pt{}=\Bigg(
-\sum_{\alpha\ne\beta,j\ne k}|\theta_{\alpha\alpha jk}-\theta_{\beta
\beta jk}|^2-(n+1)\sum_{\alpha,\beta,j,k,\,{\rm pairwise}\,\ne}
|\theta_{\alpha\beta jk}|^2\cr
&&{}-(4n+8)\sum_{\alpha\ne j<k\ne\alpha}|\theta_{\alpha\alpha jk}|^2
-4\sum_{\alpha\ne j}|\theta_{\alpha\alpha\alpha j}|^2\cr
&&{}-\sum_{\alpha\ne\beta\ne j\ne\alpha}|\theta_{\alpha\alpha jj}-
\theta_{\beta\beta jj}|^2-\sum_{\alpha\ne\beta}
|\theta_{\alpha\alpha\alpha\alpha}-2\theta_{\alpha\alpha\beta\beta}|^2
\\
&&{}+n\sum_\alpha\rho_{\alpha,\varepsilon}^2-
\bigg(\sum_\alpha\rho_{\alpha,\varepsilon}\bigg)^2~~\Bigg)~\omega_\varepsilon^n
\end{eqnarray*}
where all terms in the summation are nonpositive except the ones involving
the Ricci eigenvalues $\rho_{\alpha,\varepsilon}$. However, the assumption
$(**)^p$ for $p\geq 2$ implies that the integral of the difference
appearing in the last line converges to~$0$ as $\varepsilon\to 0$.
Theorem~\ref{stability} is proved.

\vskip12pt\noindent
Our results lead to the following interesting question.

\begin{question} Assuming that a ``sufficiently good'' family of
regularization operators~$J_\varepsilon$ is used, can one infer that the
resulting family $(\omega_\varepsilon)$ of fixed points such that
$\Ricci(\omega_\varepsilon)=J_\varepsilon(\omega_\varepsilon)$
$($or some subsequence$\,)$ is well behaved in $L^1$, resp.\
$L^2$ norm$\,?$
\end{question}

\begin{remark} {\rm A strategy to attack the Campana-Peternell conjecture
could be as follows. The first step would be to prove
that $H^0(X,T_X)\neq 0$ (if $\dim X>0$). Assume therefore 
$H^0(X,T_X)=0$ and try to reach a contradiction. It is well known in that
case that there are much fewer obstructions to the existence of
K\"ahler-Einstein metrics, for instance one has the Bando-Mabuchi uniqueness
theorem \cite{BM85} and (obviously) the vanishing of the classical
Futaki invariant \cite{Fut83};
hence we could expect our family of weakly almost K\"ahler-Einstein
metrics to converge to a genuine K\"ahler-Einstein metric. But then the
resulting known Chern class inequalities (the ones of Theorem~1.3~(ii)
combined with the Fulton-Lazarsfeld inequalities \cite{FL83}, \cite{DPS94})
might help to contradict $H^0(X,T_X)=0$, e.g.\ by the Riemann-Roch
formula or by ad hoc nonvanishing theorems such as the
generalized hard Lefschetz theorem~\cite{DPS01}. Now, the existence of
holomorphic vector fields on $X$
implies that $X$ has a non trivial group of automorphisms $H=\Aut(X)$,
and one could then try to apply induction on dimension on a suitable
desingularization $Y$ of $X\mathop{/\!/}H$, if $Y$ is not reduced to a
point and one can achieve $Y$ to have $T_Y$ nef.}
\end{remark}

\bibliographystyle{amsalpha}

\vskip3mm\noindent
Jean-Pierre Demailly\\
Universit\'e Grenoble Alpes, Institut Fourier, F-38000 Grenoble, France\\
CNRS, Institut Fourier, F-38000 Grenoble, France\\
\emph{e-mail}\/: jean-pierre.demailly@univ-grenoble-alpes.fr

\end{document}